\documentclass[a4paper]{article}
\usepackage[utf8]{inputenc}
\usepackage{pdfpages}
\usepackage{multicol}
\usepackage{tikz}
\usepackage[utf8]{inputenc}
\usepackage{pdfpages}
\usepackage{multicol}
\usepackage{tikz}
\newcounter{casenum}

\usetikzlibrary{patterns,positioning}
\usepackage{caption}
\usepackage{subcaption}
\usepackage{todonotes}
\usepackage{float}
\usepackage{amsmath}
\usepackage{amsfonts,amssymb,amsmath,amsthm}

\tikzstyle{every node}=[circle, draw, fill=black, inner sep=0pt, minimum width=4pt]

\usepackage[justification=centering]{caption}
\usetikzlibrary{patterns,positioning}

\usepackage{caption}
\usepackage{subcaption}
\usepackage{todonotes}
\usepackage{float}
\usepackage{amsmath}
\usepackage{amsfonts,amssymb,amsmath,amsthm}
\usepackage{anysize} 
\marginsize{2.5cm}{2.5cm}{2.5cm}{2.5cm}

\usetikzlibrary{shapes.geometric}
\usepackage{stmaryrd}
\def\mytransformation{
\pgfmathsetmacro{\myX}{0.9*\pgf@x}
\pgfmathsetmacro{\myY}{0.0045*(\pgf@x)*(\pgf@x)+(\pgf@y)}
\setlength{\pgf@x}{\myX pt}
\setlength{\pgf@y}{\myY pt}
}

\makeatother
\usetikzlibrary{backgrounds}
\usetikzlibrary{shapes.geometric}
\usepackage{enumerate}
\usepackage{mathrsfs}
\usepackage{amssymb}
\usepackage{graphicx}%
\usepackage{multirow}%
\usepackage{amsmath,amssymb,amsfonts}%
\usepackage{amsthm}%
\usepackage{mathrsfs}%
\usepackage[title]{appendix}%
\usepackage{xcolor}%
\usepackage{textcomp}%
\usepackage{manyfoot}%
\usepackage{booktabs}%
\usepackage{algorithm}%
\usepackage{algorithmicx}%
\usepackage{algpseudocode}%
\usepackage{listings}%
\usepackage{setspace}
\usepackage{pdflscape}
\usepackage{hyperref}
\usepackage{graphicx}
\usepackage{pst-node}
\usepackage{booktabs}
\usepackage{tikz}
\usepackage{authblk}
\usepackage{todonotes}

\renewenvironment{proof}{{\it Proof.}}{\qed \medskip}
\newtheorem{theorem}{Theorem}[section]

\newtheorem{obs}[theorem]{Observation}

\newtheorem{lemma}[theorem]{Lemma}
\newtheorem{definition}[theorem]{Definition}
\newtheorem{corollary}[theorem]{Corollary}

\usepackage{mathtools}

\providecommand{\keywords}[1]{\textbf{\textit{Keywords---}} #1}

\begin{document}

\title{Contractions in perfect graphs}


\author[1]{Alexandre Dupont-Bouillard}
\author[1]{Pierre Fouilhoux}
\author[1]{Roland Grappe}
\author[1]{Mathieu Lacroix}

\affil[1]{LIPN, CNRS UMR 7030, Université Sorbonne Paris Nord, France\footnote{
dupont-bouillard, fouilhoux, grappe, lacroix@lipn.fr}}

\maketitle


\abstract{
In this paper, we characterize the class of {\em contraction perfect} graphs which are the graphs that remain perfect after the contraction of any edge set.
We prove that a graph is contraction perfect if and only if it is perfect and the contraction of any single edge preserves its perfection.
This yields a characterization of contraction perfect graphs in terms of forbidden induced subgraphs, and a polynomial algorithm to recognize them.
We also define the utter graph $u(G)$ which is the graph whose stable sets are in bijection with the co-2-plexes of $G$, and prove that $u(G)$ is perfect if and only if $G$ is contraction perfect. 
}

\section*{Introduction}

A graph $G$ is called \emph{perfect} if $\omega(H)=\chi(H)$ for every induced subgraph $H$ of $G$, where $\omega(H)$ is the clique number of $H$ and $\chi(H)$ its chromatic number. 
Contracting any pair of vertices, that is, identifying both vertices, does not always preserve perfection. 
For instance, perfection is destroyed if there is an induced path of odd length between those vertices.

Two nonadjacent vertices in a graph form an {\em even pair} if every induced path between them has an even number of edges. 
Fonlupt and Uhry~\cite{FONLUPT198283} proved that contracting an even pair in a perfect graph preserves perfection, and Meyniel proved what is called the Even Pair Lemma~\cite{MEYNIEL1987313}: no minimally imperfect graph contains an even pair.
Concerning the associated decision problems, Bienstock~\cite{BIENSTOCK199185} proved that the following are CoNP-complete in the general case: deciding whether a given pair of vertices of a graph form an even pair, and deciding whether a given graph contains an even pair.
Nevertheless, in perfect graphs, both problems become solvable in polynomial time~\cite{DBLP:journals/combinatorica/ChudnovskyCLSV05}. A graph is perfectly contractile if for every induced subgraph there exists a sequence of even pair contraction yielding a clique. In~\cite{trot}, a subclass of perfectly contractile graph is given which leads to a polynomial combinatorial algorithm for the coloring problem for that class.

\medskip

In this paper, we approach contractions in perfect graphs with a slightly different point of view: we are interested in the graphs for which the contraction of any set of edges preserves perfection, and we call them {\em contraction perfect}.
Since perfection is preserved under taking induced subgraphs, contraction perfect graphs can be seen as the graphs in which every induced minor is perfect, where an {\em induced minor} is a graph obtained by deleting vertices and contracting edges. A class of graph is $\emph{induced minor closed}$ if, for any element of that class, each of its induced minors also belongs to that class. 

\medskip

A {\em co-$2$-plex} in a graph is a set of vertices inducing a subgraph of maximum degree at most one.
This notion, or rather its complement version, was introduced in 1978 by Seidman and Foster~\cite{k-plex} to seek communities in a graph with more freedom than when looking for cliques---a {\em clique} is a set of pairwise adjacent vertices.
Indeed, a {\em $k$-plex} is a set $W$ of vertices inducing a graph where every vertex has degree at least $|W|-k$, and cliques are special cases of $2$-plexes.
The underlying optimization problem is to find a maximum weight $k$-plex in a given weighted graph.
For any fixed $k$ and hence for $k=2$, this problem is NP-hard~\cite{bala}.
Hence, by complementing the graph, so is the problem of finding a maximum co-2-plex.
It turns out that this latter optimization problem, which can be seen as a relaxation of the maximum stable set problem---a {\em stable} set is a set of vertices inducing a subgraph of maximum degree at most 0---, behaves well in contraction perfect graphs.

This topic has been especially lively the past few years, mostly in the design of practical algorithms that find a maximum weight $k$-plex.
For this problem, there are exact algorithms~\cite{Zhou_Hu_Xiao_Fu_2021}, algorithms based on local search~\cite{DBLP:journals/heuristics/Pullan21,Chen2019CombiningRL,8794502,Jin2019EffectiveRL,Chen_Wan_Cai_Li_Chen_2020}, heuristics~\cite{doi:10.1142/S0218213019500155}, branch and bound algorithms, some of which incorporate machine learning ingredients~\cite{bala,10.1007/978-3-319-11812-3_21,DBLP:journals/corr/abs-2301-07300,DBLP:journals/dase/WangJYL17,ijcai2021p233,10.14778/3565816.3565817,DBLP:journals/corr/abs-2208-05763}, and quadratic models~\cite{Stetsyuk}.
From a more combinatorial optimization point of view, a combinatorial algorithm is devised in~\cite{suite2}, and polyhedral studies are conducted in~\cite{bala,pol}.
The problem of enumerating the $k$-plexes also received some attention~\cite{10.1145/3511808.3557444,WANG201779,Jabbour2022ADF}.
In~\cite{suite2}, a parameterized algorithm is proposed to both find the smallest set of vertices whose removal gives a graph with vertices of maximal degree $k$ and the maximum co-$(k\!+\!1)$-plex problem. The two latter problems parameterized by the size of the solution are not FPT reducible hence, no other results can be deduced for the maximum co-$k\!+\!1$-plex problem.

As we shall see, our approach of contractions in perfect graphs yields a family of graphs in which the co-2-plex problem is solvable in polynomial time.

\paragraph{Contributions.}
We start with a surprising result: contraction perfect graphs are characterized, among perfect graphs, as those for which the contraction of any single edge preserves perfection.
Building upon this result, we characterize contraction perfect graphs in terms of forbidden induced subgraphs.
As another byproduct, we get that the associated recognition problem, that is, the problem of deciding whether a given graph is contraction perfect, is solvable in polynomial time.

It turns out that contraction perfect graphs are graphs in which the maximum weight co-2-plex problem is solvable in polynomial time.
To prove this, we introduce the utter graph of a graph: it is a graph in which the stable sets correspond to the co-2-plexes of the original graph.
We first give sufficient conditions for an induced minor closed graph class to both contain a graph and its utter graph.
These conditions allow us to derive another characterization of contraction perfect graphs: a graph is contraction perfect if and only if its utter graph is perfect.
Since finding a maximum weight stable set can be done in polynomial time in perfect graphs~\cite{polyperfect}, this yields the announced complexity result for the maximum weight co-$2$-plex. 
We also strengthen the connection between a graph and its utter graph for split, trivially perfect, and chordal graphs.

\paragraph{Outline.} 
In Section~\ref{s1}, we first provide a forbidden induced subgraph characterization for the graphs remaining perfect after the contraction of any single edge.
We show that these graphs are exactly contraction perfect graphs in Section~\ref{ss11}.
This yields the forbidden induced subgraph characterization of contraction perfect graphs in Section~\ref{ss12}.
Consequently, deciding if a given graph is contraction perfect can be done in polynomial time.

In Section~\ref{s2}, we introduce the utter graph of a graph.
This yields another characterization of contraction perfect graphs in Section~\ref{ss21} as the graphs whose utter graphs are perfect.
As a byproduct, we get that a maximum weight co-2-plex can be found in polynomial time in a contraction perfect graphs.
In Section~\ref{ss22}, we strengthen the link with the utter graph for a few subclasses of contraction perfect graphs.

\paragraph{Definitions.}
All the graphs in this paper are simple and connected. 
Given a graph $G=(V,E)$, we denote its \textit{complement} by $\overline{G}=(V,\overline E)$, where $\overline{E} = \{ uv: uv \notin E \}$. 
We denote by $V(G)$ (resp. $E(G)$) the vertex (resp. edge) set of $G$. Two vertices $u$ and $v$ are \emph{adjacent} if $uv \in E(G)$. 
A vertex is \textit{universal} if it is adjacent to all the other vertices. Given a subset of vertices $W \subseteq V$, let $E(W)$ denote the set of edges of $G$ having both endpoints in $W$ and $\delta (W)$ the set of all edges having exactly one endpoint in $W$. When $W$ is a singleton $\{w\}$, we will simply write $\delta (w)$. We say that the edges in $\delta (w)$ are \textit{incident} to $w$, and two edges sharing an extremity are also said {\em incident}.
A \textit{matching} is a set of pairwise nonincident edges. 
For $F \subseteq E$, let $V(F)$ denote the set of vertices incident to any edge of $F$.
Given $W\subseteq V$, the graph $G[W]=(W,E(W))$ is the \textit{subgraph induced by $W$ in $G$}. When $H$ is an induced subgraph of $G$, we say that $G$ {\em contains} $H$. 
Given a vertex $u\in V$, we denote  by $N_u = \{w \in V : uw \in E \}$ its \textit{neighborhood} in $G$, and by $\overline{N}_u = N_u \cup \{ u\}$ its \textit{closed neighborhood}. Two vertices $u$ and $v$ are true (resp. false) twins if $\overline{N}_u = \overline{N}_v$ (resp. $N_u = N_v)$.

A \textit{path} (resp. \textit{hole}) is a graph induced by a sequence of vertices $(v_1,\dots,v_p)$ whose edge set is $\{ v_iv_{i+1}:  \ i = 1,\dots, p-1 \} $ (resp. $\{ v_iv_{i+1}:  \ i = 1,\dots, p-1 \} \cup \{v_1v_p \}$ with $p\ge4$). Note that this definition usually corresponds to induced paths.
A subset of vertices induces a path (resp. hole) if its elements can be ordered into a sequence inducing a path (resp. hole).
An \emph{antipath} (resp. \textit{antihole}) of $G$ is a path (resp. hole) of $\overline{G}$.
The size of a graph $H$ is $|V(H)|$, and its parity is the parity of $|V(H)|$.

The \textit{contraction} of an edge $uv$ in $G$ consists in deleting $u$ and $v$, and adding a new vertex $w$ and the edges $wz$ for all $z\in N_u \cup N_v$. 
This new graph is denoted $G/uv$.
For $F \subseteq E$, we denote by $G/F$ the graph obtained from $G$ by contracting all the edges in $F$. The \emph{image} of a vertex $v$ of $G$ in $G/F$ is the vertex of $G/F$ to which $v$ is contracted, and the \emph{image} of a set of edges $L$ is the set of the images of the vertices of $V(L)$ in $G/F$.
An edge $uv$ and a vertex $w$ are \textit{adjacent by contraction} if $w$ is adjacent to $x_{uv}$ in $G/uv$, where $x_{uv}$ is the image of $uv$ in $G/uv$. 
In other words, at least one of $uw$ and $vw$ is in $E$.
Two edges $uv$ and $xy$ are \textit{adjacent by contraction} if contracting both edges results in two adjacent vertices, that is, $\delta(\{u,v\}) \cap \delta (\{x,y \})\neq \emptyset$. 

Given a total order on a finite set of elements, an \emph{interval} is a subset of consecutive elements following that order.
\emph{The interval graph} $G_{\mathcal{I}}$ of a given finite set of intervals $\mathcal{I}$ is the graph having one vertex per interval in $\mathcal{I}$ and an edge between two vertices if their corresponding intervals have a nonempty intersection.



\section{Contraction perfect graphs}\label{s1}

This section investigates how contracting edges in a perfect graph impacts its perfection.
It relies on  the strong perfect graph theorem~\cite{perfectgraphtheorem}, which states that perfect graphs are the graphs that contain neither odd holes nor odd antiholes.
We start by characterizing when the contraction of an edge destroys the perfection of a graph.
In the second and third subsections, we provide two characterizations of contraction perfect graphs: the first one compares contractions of single edges with contractions of sets of edges, and the second one is in terms of forbidden induced subgraphs.

\medskip

We start by proving a parity result about interval graphs that will be used in Section~\ref{sec:1.1} but might be of independent interest.
A set of intervals is an \emph{odd intersection interval set} if the non-empty intersection of any subset of intervals is an interval of odd cardinality.

\medskip
\begin{lemma}\label{lem:patates}
Given an odd intersection interval set $\mathcal{I}$, the union of the intervals associated with every connected component of $G_{\mathcal{I}}$ has odd cardinality.   
\end{lemma}
\noindent\begin{proof}
It is enough to prove the result when $G_{\mathcal{I}}$ is connected.
We proceed by induction on $|\mathcal{I}|$, the case $|\mathcal{I}|=1$ being immediate.
Let $I$ and $I'$ be two intersecting intervals of $\mathcal{I}$. Hence, $I\cup I'$ is an interval. Let $\mathcal{I'}$ be obtained from $\mathcal{I}$ by replacing $I$ and $I'$ by $I \cup I'$. Since $\mathcal{I}$ is an odd intersection interval set, $|(I \cup I') \cap K| = |I \cap K| + |I' \cap K| - |I \cap I' \cap K|$ is odd or zero for each intersection $K$ of any subset of intervals of $\mathcal{I}$. Hence, $\mathcal{I}'$ is also an odd intersection interval set. Moreover, $G_{\mathcal{I'}}$ is connected as it is obtained from $G_{\mathcal{I}}$ by contracting the edge whose extremities are the vertices associated with $I$ and $I'
$.
Since $|\mathcal{I}'| < |\mathcal{I}|$, the induction hypothesis implies that the union of intervals of $\mathcal{I}'$ has odd cardinality. By definition, the latter set is equal to the union of the intervals of $\mathcal{I}$, which ends the proof.
\end{proof}


\subsection{When does the contraction of an edge destroy perfection? \label{sec:1.1}}

By the strong perfect graph theorem, if the contraction of an edge $e$ in a perfect graph $G$ destroys its perfection, this implies that $G/e$ contains an odd hole or an odd antihole. Moreover, since $G$ is perfect, the image of $e$ must be a vertex of such a forbidden induced subgraph. Lemmas~\ref{lem:holegt5} and~\ref{lem:expandedAntiholes} characterize when contracting an edge yields an odd hole or an odd antihole, respectively. Both results will be used in the proof of Theorem~\ref{thm:holeantihole} which characterizes when the contraction of an edge destroys the perfection of a graph.

\medskip

\begin{lemma}\label{lem:holegt5}
If $G/F$ contains a hole $H$ for some edge set $F$, then $G$ contains a hole of size at least~$|V(H)|$.
\end{lemma}


\noindent\begin{proof}
We proceed by induction on $|F|$.
The result holds if $|F|=0$.
Otherwise, some vertex $v$ of $H$ is obtained by contracting the edges of a connected subgraph $(V',F')$ of $G$ with $\emptyset\neq F'\subseteq F$.
Let $a$ and $b$ be the neighbors of $V'$ in $H$, and let $P$ be an $ab$-path in $G[V'\cup\{a,b\}]$.
Then, $V(H)\setminus \{v\}\cup V(P)$ induces a hole of $G/(F\setminus F')$ of size at least $|V(H)|$, and induction concludes.
\end{proof}

As we shall see in Lemma~\ref{lem:expandedAntiholes}, when the contraction of an edge destroys the perfection of $G$ by creating an odd antihole, this implies that $G$ contains the following specific structure.

\medskip

\begin{figure}[h]
\centering

	\begin{tikzpicture}[state/.style={circle, minimum size=0.8cm},scale = 1]
 
            \draw [black,,fill=gray!30] plot [smooth cycle] coordinates {
            (-4.5,3) (-4.7,4.2)  (-3.5,4.4)  (0,1.8)  
            (3.5,4.4) (4.7,4.2)  (4.5,3) (0,0)
            };

		      \node[state,draw,circle,fill=gray!70] (1)at(-1,6) {$u$};
            \node[state,draw,circle,fill=gray!70] (2)at(1,6) {$v$};
            \node[state,draw,circle,fill=gray!70] (3)at(4,3.7) {$w_{p}$};
            \node[state,draw,circle,fill=gray!70,scale = 0.83] (4)at(3,2.5) {$w_{p-1}$};
            \node[state,draw,circle,fill=gray!70,] (5)at(-3,2.5)  {$w_{2}$};
            \node[state,draw,circle,fill=gray!70] (6)at(-4,3.7)  {$w_1$};
             \node[circle,draw,fill=gray!70,ellipse, minimum width = 1.4cm,minimum height = 0.8cm,text opacity = 1] (7) at(0,1) {$ w_{3},\dots,w_{p-2}$};
             
            \draw (1)-- (2);
            \draw (2)--(4);
            \draw (1)-- (5);
            \draw  [line width = 0.5mm](1) -- (7);
            \draw [line width = 0.5mm](2) -- (7);

	\end{tikzpicture}

\caption{An expanded antihole.}
\label{Loah}
\end{figure}
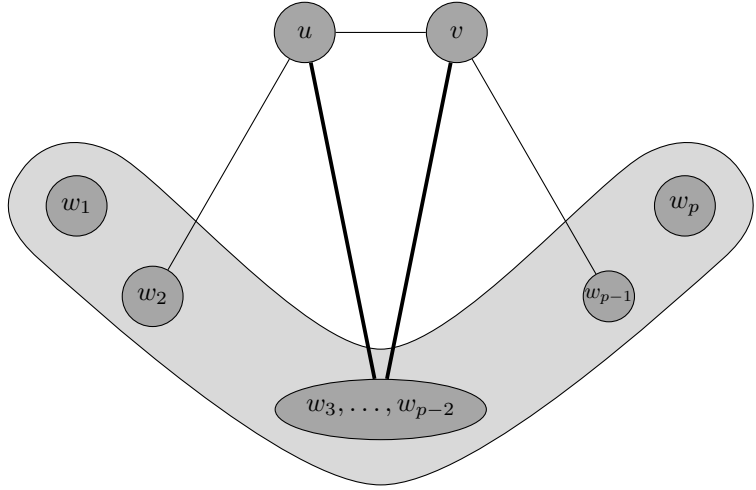

\begin{definition}\label{def:eah}
An edge $e$ and an even antipath $P$ induced by $(w_1,\dots,w_p)$ with $p\ge 6$ form an \emph{expanded antihole} if one extremity of $e$ is adjacent to all vertices of $P$ but $w_1, w_{p-1}, w_p$, and the other extremity is adjacent to all vertices of $P$ but $w_1, w_2, w_p$. 
\end{definition}
\medskip
Figure~\ref{Loah} represents an expanded antihole where $(w_1,\dots,w_p)$ induces an even antipath. An edge $e$ is \emph{involved} in an expanded antihole if there exists an even antipath forming with $e$ an expanded antihole. 
Note that, by definition of an expanded antihole, the extremities of $e$ are not vertices of $P$, and contracting $uv$ yields an odd antihole.
Moreover, if $e=uv$, either $(u,w_2,\dots,w_{p-1},v)$ or $(v,w_2,\dots,w_{p-1},u)$ induces an even antipath that forms an expanded antihole with edge $w_1w_p$.

            
             

Note that if one relaxes the adjacency requirement in the definition an expanded antihole to: "each vertex $w_2,\dots,w_{p-1}$ is adjacent to one extremity of $e$", contracting $e$ still yields an odd antihole of $G/e$. However, in this case, $G$ is not necessarily perfect.
The following lemma asserts that the only case where $G$ is perfect and $G/e$ contains an odd antihole of size at least 7 is when $G$ contains an expanded antihole.

\begin{lemma}\label{lem:expandedAntiholes}
Let $G$ be a perfect graph, $H$ a subgraph of $G$, and $uv$ an edge of $E(H)$. $H/uv$ is an odd antihole of $G/uv$ of size at least 7 if and only if $H$ is an expanded antihole formed by $uv$ and the antipath induced by $V(H)\setminus \{u,v \}$.
\end{lemma}    
\noindent\begin{proof}
$(\Leftarrow)$
    By definition, if $H$ is an expanded antihole involving $uv$, then $H/uv$ is an odd antihole of $G/uv$.
    
\noindent$(\Rightarrow)$
Let us denote by $w_{uv}$ the vertex obtained by contracting $uv$. Since $H/uv$ is an antihole, $V(H)\setminus \{u,v\}$ can be ordered into the sequence $(w_1,\dots,w_{p})$ with $p\ge 6$ inducing an even antipath in $G$. The set of vertices of a maximal antipath of $H[N_{u} \setminus v  ]$ (resp. $H[N_v \setminus u]$) forms an interval of $(w_2,\dots,w_{p-1})$. Let $\mathcal{P}^u$ (resp. $\mathcal{P}^v$) be the set of these intervals, and $\mathcal{I} = \mathcal{P}^u \cup \mathcal{P}^v$. Since $G$ is perfect, $H$  contains no antihole and then, the interval $( w_2,\dots , w_{p-1})$ does not belong to $\mathcal{I}$.

Suppose that interval graph $G_{\mathcal{I}}$ associated to $\mathcal{I}$ is not connected. Then, there exists $i \in \{2,\dots, p-2 \}$ such that no interval of $\mathcal{I}$ contains both $w_i$ and $w_{i+1}$. By definition of $\mathcal{I}$, $\{u,v,w_i, w_{i+1} \}$ induces a path of $H$ with extremities $w_i$ and $w_{i+1}$. Since $p \ge 6$, there exists $z \in \{ w_1,w_p\}$ adjacent to both $w_i$ and $w_{i+1}$. As $z$ is adjacent to neither $u$ nor $v$, $\{z,w_i,u,v,w_{i+1} \}$ induces an odd hole of $G$, a contradiction to its perfection. Hence, $G_{\mathcal{I}}$ is connected.

Note that every interval $I = \{ w_i, \dots, w_j\}$ of $\mathcal{I}$ has odd cardinality, as otherwise $\{ w_{i-1},w_i, \dots, w_j,$ $ w_{j+1}\}$ would induce, together with $u$ if $I \in \mathcal{P}^u$ or $v$ if $I \in \mathcal{P}^v$, an odd antihole of $G$, contradicting its perfection. 
By Lemma~\ref{lem:patates}, since $p$ is even and $G_{\mathcal{I}}$ is connected, $\mathcal{I}$ is not an odd intersection interval set. Moreover, the intersection of more than two intervals of $\mathcal{I}$ is empty by definition of $\mathcal{P}^u$ and $\mathcal{P}^v$. 
Hence, there exists $I \in \mathcal{P}^u$ and $I' \in \mathcal{P}^v$ whose intersection is nonempty and has even cardinality. Let $I = \{w_i,\dots,w_k \}$ and $I' = \{w_j,\dots,w_\ell \}$ with $2 \le i \le j \le k \le \ell \le p-1$, Figure~\ref{proof_theorem} gives an illustration. Since $|I|$ and $|I'|$ are odd and $|I \cap I'|$ is even, $2 \le i < j < k <\ell \le p-1$. Suppose that $j > 3$. Then, $w_1$ is adjacent to $w_{j-1}, \dots, w_{k+1}$ so $\{w_1, v,w_{j-1}, \dots, w_{k+1},u \}$ induces an odd antihole of $G$, a contradiction. Hence, $j=3$. Similarly, one can prove that $k=p-2$. Hence, $u$ is adjacent to $w_2,\dots, w_{p-2}$ and $v$ is adjacent to $w_3, \dots, w_{p-1}$. Therefore, $(w_1,\dots,w_p)$ and the edge $uv$ form an expanded antihole.
\end{proof}

\begin{center}
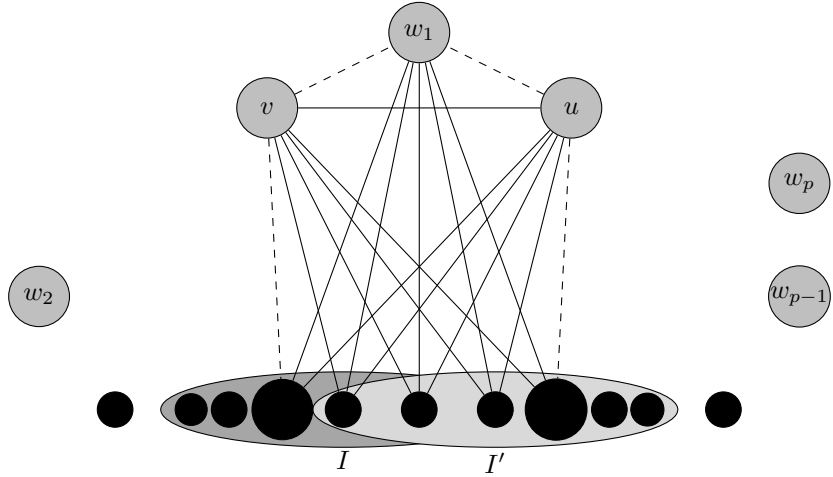
\begin{figure}[h]
\centering
	\begin{tikzpicture}[state/.style={circle, minimum size=0.8cm},scale = 1]
	      \node[state,draw,circle,fill=gray!50]        (1)at(-2,7) {$v$};
            \node[state,draw,circle,fill=gray!50] (2)at(2,7) {$u$};
            \node[state,draw,circle,fill=gray!50] (3)at(0,8) {$w_{1}$};
            \node[state,draw,circle,fill=gray!50,scale = 1] (4)at(5,4.5) {$w_{p-1}$};
            \node[state,draw,circle,fill=gray!50] (5)at(-5,4.5)  {$w_{2}$};
            \node[state,draw,circle,fill=gray!50] (6)at(5,6)  {$w_p$};
             \node[circle,draw,fill=gray!70,ellipse, minimum width = 4.8cm,minimum height = 1cm,text opacity = 1,label = below:$I$] (7) at(-1,3) {};
			  \node[circle,draw,fill=gray!30,ellipse, minimum width =4.8cm,minimum height = 1cm,text opacity =1,opacity = 0.5,,label = below:$I'$] (8) at(1,3) {};
				\node[text opacity = 1, opacity= 0] (9) at (3,3) {$w_{\ell}$};
                \node[text opacity = 1, opacity= 0] (18) at (4,3) {$\dots$};
                \node[text opacity = 1, opacity= 0] (10) at (2.5,3) {$\dots$};
                \node[text opacity = 1, opacity= 0] (11) at (1.8,3) {$w_{k+1}$};
                \node[text opacity = 1, opacity= 0] (12) at (1,3) {$w_{k}$};
                \node[text opacity = 1, opacity= 0] (13) at (0,3) {$\dots$};
                \node[text opacity = 1, opacity= 0] (14) at (-1,3) {$w_{j}$};
                \node[text opacity = 1, opacity= 0] (15) at (-1.8,3) {$w_{j-1}$};
                \node[text opacity = 1, opacity= 0] (16) at (-2.5,3) {$\dots$};
                \node[text opacity = 1, opacity= 0] (19) at (-4,3) {$\dots$};
              \node[text opacity = 1, opacity= 0] (17) at (-3,3) {$w_{i}$};
              \draw [dashed] (1) -- (15);
              \draw  (1) -- (14);
              \draw  (1) -- (13);
              \draw  (1) -- (12);
              \draw  (1) -- (11);
              \draw (3) -- (15);
              \draw [dashed] (1) -- (3);
              \draw [dashed] (2) -- (3);
            \draw (1)-- (2);
            \draw (3) -- (11);
            \draw (3) -- (12);
            \draw (3) -- (13);
            \draw (3) -- (14);
            \draw [dashed](2) -- (11);
            \draw (2) -- (12);
            \draw (2) -- (13);
            \draw (2) -- (14);
            \draw (2) -- (15);
	\end{tikzpicture}
	
	\caption{Illustration of the proof of Lemma~\ref{lem:expandedAntiholes} (the doted lines \\ represent nonedges and dots the remaining vertices of the antipath).}
	\label{proof_theorem}
	\end{figure}
\end{center}
\begin{theorem}\label{thm:holeantihole}
Let $G=(V,E)$ be a perfect graph and $uv \in E$. The graph $G/uv$ is not perfect if and only if $G$ contains an even hole of size at least 6 containing $uv$ or an expanded antihole involving $uv$. 
\end{theorem}
\noindent\begin{proof}
$(\Leftarrow)$
In both cases, as mentioned before, contracting $uv$ destroys perfection.

\noindent$(\Rightarrow)$
Suppose that $G/uv$ is not perfect.
If $G/uv$ contains an odd hole $H$, then $G$ contains a hole $H'$ of size at least $|V(H)|$ by Lemma~\ref{lem:holegt5}. Since $G$ is perfect, $|V(H')|$ is even and at least 6, and contains $uv$. Otherwise, $G/uv$ contains an odd antihole hence, by Lemma~\ref{lem:expandedAntiholes}, $G$ contains an expanded antihole involving $uv$.
\end{proof}


\subsection{Contracting an edge or a set of edges?}\label{ss11}

We prove that contraction perfect graphs are characterized by the contraction of single edges.
\begin{theorem} \label{the:contractequiv}
A perfect graph is contraction perfect if and only if it remains perfect by the contraction of any single edge.
\end{theorem}
\noindent\begin{proof}
To prove the nontrivial direction, let $G=(V,E)$ be a perfect graph and $F \subseteq E$ be such that $G/F$ is not perfect. We will prove that $G/e$ is not perfect for some $e \in E$. Without loss of generality, suppose that $F$ has minimum size, that is, $G/F'$ is perfect for all $F' \subsetneq E$ such that $|F'| < |F|$.     
Since $G/F$ is not perfect, it contains an odd hole or an odd antihole.
Suppose that $G/F$ contains an odd hole $C$. By Lemma~\ref{lem:holegt5}, $G$ contains a hole $D$ with $|V(D)|\ge |V(C)|$. Since $G$ is perfect, $|V(D)|$ is even and at least $6$. Then, $G/e$ contains the odd hole $D/e$ for any $e \in E(D)$ and thus is not perfect.

From now on, we may assume that $G$ and $G/F'$ for all $F' \subseteq E$ contain no hole of size at least 5. Hence, $G/F$ contains odd antiholes and let us denote by $(w_0,\dots,w_p)$ one of them where $p\ge 6$.

%
We will first prove that $F$ is a matching. To do this, let us suppose that $F$ contains two adjacent edges $uv$ and $vw$ whose image in $G/F$ is without loss of generality $w_0$. Let $G' = G/(F \setminus \{ uv,vw\})$. Since $G'/\{ uv,vw\}$ contains an odd antihole, and since $G'/uv$ is perfect by minimality of $F$ and has no hole of size at least 5, then, by Theorem~\ref{thm:holeantihole}, $G'/uv$ contains an expanded antihole formed by the even antipath $P_0 = (w_1,\dots w_p )$, and the edge $x_{uv}w$ where $x_{uv}$ is the image of $uv$ in $G'$. 

Without loss of generality suppose that $x_{uv}$ is adjacent to $\{w_2,\dots,w_{p-2} \}$ and $w$ is adjacent to $\{ w_3,\dots,w_{p-1}\}$. Note that $v$ is not adjacent to $w_2$ in $G$, as otherwise $G'/vw$ contains the odd antihole induced by $\{x_{vw},w_1,\dots, w_p \}$ ---where $x_{vw}$ is the image $vw$ in $G'/vw$--- which contradicts the minimality assumption on $F$.
Now, $G'$ contains $uw$ as otherwise $\{u,v,w,w_{p-1},w_2 \}$ induces an odd hole of $G'$, contradicting its perfection. But now, $G'/uw$ contains the odd antihole induced by $\{ x_{uw}, w_1,\dots,w_p\}$ where $x_{uw}$ is the vertex obtained by contracting edge $uw$. This contradicts the minimality assumption on $F$, hence $F$ is a matching.

We now prove that $F$ only contains one edge. Let $uv$ and $u'v'$ be two edges of $F$. By minimality assumption, their images are vertices of $\{w_0,\dots,w_p\}$, and up to node relabelling, let us suppose that the image of $uv$ is $w_0$ and the one of $u'v'$ is $w_i$ with $i \in \{1,\dots,p/2 \}$. Set $G' = G/(F\setminus \{uv, u'v' \})$. By Lemma~\ref{lem:expandedAntiholes}, $uv$ forms with $P_0 = \{w_1,\dots,w_p \}$ an expanded antihole of $G'/u'v'$. Similarly, $u'v'$ forms with $P_i =\{w_{i+1},\dots,w_p,w_0,\dots,w_{i-1}  \}$ an expanded antihole of $G'/uv$. By definition of expanded antiholes, we suppose without loss of generality that $u$ (resp. $u'$) is adjacent to $w_2$ (resp. $w_{i-2}$) but not to $w_{p-1}$ (resp. $w_{i+2}$). 
Hence, $H_u = \{u,w_1,\dots,w_{p-1} \}$ induces an even antihole of $G'/u'v'$, similarly $\{w_0,\dots,w_{i-1}, u',w_{i+2},\dots,w_p, \}$ induces an even antihole of $G'/uv$. Figure~\ref{fig:oddantihole} gives an illustration of $H_u \cup u' \setminus \{w_i, w_{i+1}\}$, which induces an antihole in $G'$ of size $p-1$, and hence is odd. This contradicts the minimality assumption on $F$. Therefore $F$ cannot contain two edges, so $G$ does not remain perfect by the contraction of the edge of $F$.
  \end{proof}

  \tikzset{thick/.style={line width=1.2pt}}
		\begin{figure}
  \centering
					 \begin{tikzpicture}[state/.style={circle, draw, minimum size=0.8cm},scale = 0.3]
            \node[ellipse, label=above:{$w_0$}, fill=white,minimum width = 2cm,minimum height = 1cm] (13) at (6,12) {};
            \node[state,draw,circle,fill=gray!70,scale = 0.7] (1)at(8,12) {$v$};
            \node[state,draw,circle,fill=gray!30,scale = 0.7] (2)at(4,12) {$u$};
            \node[state,draw,circle,fill=gray!30,scale = 0.7] (3)at(0,10) {$w_1$};
            \node[scale = 1.2,fill=gray!30,text opacity=1,rotate=83,ellipse,minimum width = 1.3cm,minimum height = 0.4cm] (11)at(-2,6)  {$\cdots$};
            \node[state,draw,circle,fill=gray!30,scale = 0.7] (5)at(-2,1) {$w_{i-1}$};
            \node[ellipse, label=above:{$w_i$},rotate=152, fill=white,minimum width = 1.7cm,minimum height = 1cm] (10) at (1.5,-3.25) {};
            \node[state,draw,circle,fill=gray!30,scale = 0.7] (6)at(0,-2.5)  {$u'$};
            \node[state,draw,circle,fill=gray!70,scale = 0.7] (7)at(3,-4)  {$v'$};
            \node[state,draw,circle,fill=gray!70,scale = 0.7] (8)at(7.5,-4.5)  {$w_{i+1}$};
            \node[state,draw,circle,fill=gray!30,scale = 0.7] (9)at(11.5,-3.5)  {$w_{i+2}$};
            \node[scale = 1.2,fill=gray!30,text opacity=1,rotate=80,ellipse,minimum width = 1.3cm,minimum height = 0.4cm] (10)at(14,2)  {$\cdots$};
            \node[state,draw,circle,fill=gray!30,scale = 0.7] (11)at(14,7)  {$w_{p-1}$};
            \node[state,draw,circle,fill=gray!70,scale = 0.7] (12)at(12,10)  {$w_p$};

             \color{gray}
            \draw (2) -- (5);
            \draw (2) -- (6);
            \draw [dashed](2) -- (3);
            \draw (2) -- (9);
            \draw [dashed] (2) -- (11);
            \draw (3)-- (5);
            \draw (3) -- (6);
            \draw (3) -- (9);
            \draw (3) -- (11);
            \draw (5) -- (9);
            \draw (5) -- (11);
            \draw [dashed ](5) -- (6);
            \draw [dashed] (6) -- (9);
            \draw (9) -- (11);
            \draw (6) -- (11);

        \end{tikzpicture}
	
	\caption{Illustration for the proof of Theorem~\ref{the:contractequiv}: \\ the odd antihole $H_u \cup u' \setminus \{w_i, w_{i+1}\}$ of $G'$ is given in light gray }
 \label{fig:oddantihole}

 \end{figure}
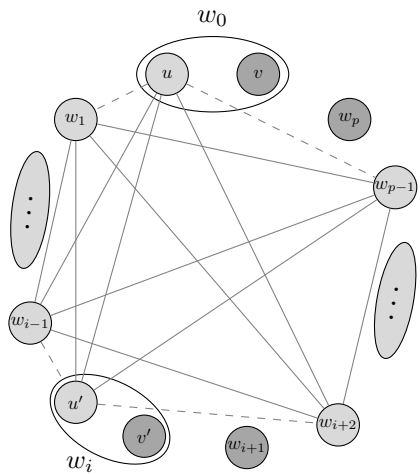

\subsection{Characterization by forbidden induced subgraphs}\label{ss12}

By definition, a perfect graph remains perfect by any vertex set deletion.
Thus, a forbidden induced subgraph characterization may be hoped for, and this tremendous result is the strong perfect graph theorem. 
Since perfect graphs do not remain perfect by any edge set contraction (a hole of size 6---which is perfect---does not remain perfect by the contraction of a single edge), there is no characterization of perfect graphs by forbidden induced minors.

By definition of contraction perfect graphs, a characterization in terms of forbidden induced minors can be directly deduced from the strong perfect graph theorem. Indeed, if a graph is not contraction perfect, then there exists an edge set contraction and a vertex set deletion that gives an odd hole or an odd antihole.
Therefore, odd hole and odd antiholes are the induced minors to be forbidden in order to ensure contraction perfection. This list of forbidden induced minors is inclusionwise minimal because of the strong perfect graph theorem.

\medskip

A list of forbidden induced subgraphs that characterizes contraction perfect graphs would be the set of graphs that yield an odd hole or an odd antihole by an edge set contraction.
Using Theorem~\ref{the:contractequiv}, which limits the set of such candidates, this section gives the minimal forbidden induced subgraph characterization of contraction perfect graphs.

\medskip

\begin{corollary}\label{cor:characCP}
A graph is contraction perfect if and only if it contains no hole of size at least 5, no odd antihole, and no expanded antihole.
\end{corollary}
\noindent\begin{proof}
$(\Rightarrow)$ By Theorem~\ref{thm:holeantihole} and the strong perfect graph theorem, a graph containing one of these structures is not contraction perfect. 

\noindent$(\Leftarrow)$ If the graph is not perfect, then it contains an odd hole or an odd antihole. If the graph is perfect but not contraction perfect, then there exists an edge whose contraction destroys the perfection of the graph by Theorem~\ref{the:contractequiv}. By Theorem~\ref{thm:holeantihole}, this implies that the graph contains an even hole of size at least 6 or an expanded antihole.
\end{proof}

The next lemma gives an example of a non trivial contraction perfect graph class and will be used as a technical lemma. Remark that $H$ is an induced subgraph of $G$ if and only if $\overline{H}$ is an induced subgraph of $\overline{G}$.

\begin{lemma}\label{lem:antipath}
    Antipaths are contraction perfect.
\end{lemma}
\noindent\begin{proof}
    The complement $\overline{P} $ of an antipath $P$ is a path. Since the path $\overline{P}$ is perfect, $P$ is perfect by the weak perfect graph theorem~\cite{lo72}. An antipath only contains holes of size at most 4.  
    An antipath contains no antihole since $\overline{P}$ contains no hole. 
    Moreover, an antipath contains no expanded antihole. Indeed, 
    any contraction of an edge $uv$ in $P$ yields a vertex whose neighborhood has size at least $|V(P)|-2$. Hence the degree of such a vertex is too large for the latter to belong to an antihole. Hence, by Lemma~\ref{lem:expandedAntiholes}, $P$ does not contain any expanded antihole involving $uv$. The result holds by Corollary~\ref{cor:characCP}.
\end{proof}

Note that that the list of forbidden induced subgraphs of Corollary~\ref{cor:characCP} is inclusionwise minimal. Indeed, if it were not the case, then one forbidden induced subgraph of Corollary~\ref{cor:characCP} would contain another one as a proper induced subgraph, which is not the case as shown in the next observation. 

A graph is {\em minimally non contraction perfect} if it is not contraction perfect and each of its proper induced subgraphs is contraction perfect.

\begin{obs}
Expanded antiholes, odd antiholes and holes of size at least 5 are minimally non contraction perfect.
\end{obs}
\noindent\begin{proof}
Every proper induced subgraph of a hole is a set of disjoint paths and hence is contraction perfect. Every proper induced subgraph of an odd antihole is a set of fully connected antipaths, the complement of such a graph is a set of disjoint paths that is contraction perfect.

We will show that expanded antiholes are minimally non contraction perfect by considering the complement of an expanded antihole.  
Let us consider $G$ an expanded antihole. Note that, $\overline{G}$ contains no antihole of size more than~6 since it contains precisely~6 vertices whose degree is greater or equal to~3, that are $W =\{u,v,w_1,w_2,w_{p-1},w_{p} \}$. Indeed, since an antihole of size at least 6 is a regular graph whose vertices have degree at least 3, these are the only candidates that may belong to such an antihole.
 Since $\overline{G}[W]$ is not regular, it is not an antihole of size 6, this proves that $G$ does not contain any holes of size at least 6.
Obviously, $\overline{G}$ does not contain any complement of an expanded antihole as proper induced subgraphs.
Now, note that the holes of $\overline{G}$ are even, which means that $G$ only contains even antiholes and hence, no hole of size 5. By Corollary~\ref{cor:characCP}, expanded antiholes are contraction perfect.
\end{proof}

\begin{corollary}\label{recognizing}
Recognizing contraction perfect graphs can be done in polynomial time.
\end{corollary}
\noindent\begin{proof}
By Theorem~\ref{the:contractequiv}, deciding whether a graph $G=(V,E)$ is contraction perfect amounts to check if $G$ is perfect and if, for each contraction of a single edge, the resulting graph is perfect. 
Each of these $|E|+1$ perfection tests can be done in polynomial time~\cite{polyperfect}. 
\end{proof}


\section{Utter graph and co-2-plexes}\label{s2}


\subsection{Perfection of the utter graph}\label{ss21}

Let $S$ be a co-2-plex. The vertex edge representation of a $S$ is a couple $(W,F)$ where $W$ are the isolated vertices of $G[S]$ and $F$ are its isolated edges. 
 Then, by definition of co-2-plexes, 
 $F$ is a matching and contracting $F$ gives $|F|$ isolated vertices nonadjacent to $W$. 
 In other words, $W$ and the image of $F$ in $G/F$ form a stable set of $G/F$ of size $|W|+ |F|$. 

\medskip

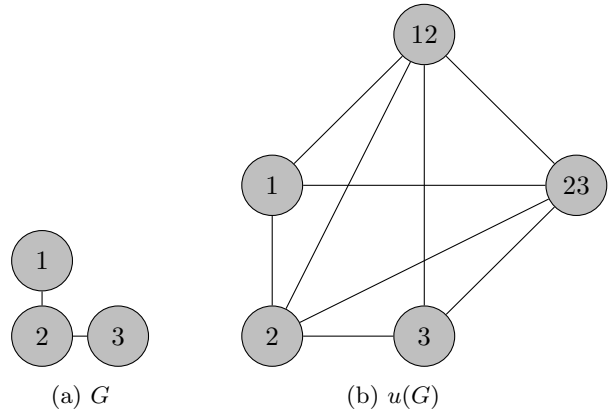
\begin{figure}[!ht]
\hspace{3.5cm}
\begin{subfigure}[b]{0.25\textwidth}
	\center
\begin{tikzpicture}[state/.style={circle, draw, minimum size=0.8cm},scale = 0.50]
            \node[state,draw,circle,fill=gray!50] (1)at(0,2) {$1$};
            \node[state,draw,circle,fill=gray!50] (2)at(0,0) {$2$};
            \node[state,draw,circle,fill=gray!50] (3)at(2,0) {$3$};
            \draw (1) -- (2);
			\draw (2) -- (3);
\end{tikzpicture}
\caption{$G$}
\end{subfigure}
\begin{subfigure}[b]{0.25\textwidth}
	\center
\begin{tikzpicture}[state/.style={circle, draw, minimum size=0.8cm}]
            \node[state,draw,circle,fill=gray!50] (1)at(0,2) {$1$};
            \node[state,draw,circle,fill=gray!50] (2)at(0,0) {$2$};
            \node[state,draw,circle,fill=gray!50] (3)at(2,0) {$3$};
            \node[state,draw,circle,fill=gray!50,] (4)at(4,2) {$23$};
            \node[state,draw,circle,fill=gray!50] (5)at(2,4) {$12$};
            \draw (4) -- (1) ;
            \draw (5) -- (1);
            \draw (5) -- (2);
            \draw (5) -- (3);
            \draw (4) -- (3);
            \draw (4) -- (2) ;
            \draw (4) -- (5);
            \draw (1) -- (2);
			\draw (2) -- (3);
\end{tikzpicture}
\caption{$u(G)$}
\end{subfigure}
\caption{A graph $G$ and its utter graph $u(G)$.}
\label{color}
\end{figure}

We use this remark to define the utter graph of a graph $G$ in which the stable sets are in bijection with the co-2-plexes of $G$. The \emph{utter graph} $u(G)$ of a graph $G=(V,E)$ has vertex set $V\cup E$ and two vertices in~$u(G)$ are adjacent if and only if their corresponding elements in $G$ are either adjacent, incident, or adjacent by contraction in $G$. Figure~\ref{color} gives an illustration of this definition where the vertices $12$ and $23$ of $u(G)$ respectively correspond to the edges $12$ and $23$ of $G$.
For each edge $uv$ of $G$, $G/uv$ is the subgraph of $u(G)$ induced by $V \setminus \{ u,v \} \cup uv$ where $uv$ denotes the vertex of $u(G)$ associated with edges $uv$.
Similarly, for a matching $F \subseteq E$, $G/F$ is the subgraph of $u(G)$ induced by $V \setminus V(F)$ and the vertices of $u(G)$ associated with edges in $F$.

\medskip

\begin{lemma}\label{bij}
There is a bijection between the co-2-plexes of $G$ and the stable sets of $u(G)$.
\end{lemma}

\noindent\begin{proof}
By definition of utter graphs, a vertex subset $W$ and an edge set $F$ is a vertex edge representation $(W,F)$ if and only if $W \cup T$ is a stable set of $u(G)$.
\end{proof}

Note that by definition of contraction perfect graphs, adding a true twin to any vertex preserves contraction perfection. Moreover, the forbidden induced subgraph characterization also proves that adding false twins preserves contraction perfection. Indeed, none of these forbidden induced subgraphs contain false twins. This gives the following replication lemma for contraction perfect graphs. 

\begin{lemma}\label{lem:replication}
    Adding true or false twins preserves contraction perfection.
\end{lemma}

We say that a graph class is \emph{replicable} if adding true twins to any graph of that class yields another graph of that class.

\begin{theorem}\label{the:utterstable}
    Let $\mathcal{C}$ be a induced minor closed and replicable graph class, then $G \in \mathcal{C}$ if and only if $u(G) \in \mathcal{C}$.
\end{theorem}
\noindent\begin{proof}
 ($\Leftarrow$) Since $G$ is an induced subgraph of $u(G)$, the result follows from the assumptions on $\mathcal{C}$.
 
\noindent($\Rightarrow$)  
Let $\widetilde{G}$ be obtained from $G$ by adding, for every edge $uv \in E$, a twin $u'$ (resp. $v'$) to $u$ (resp. $v$). By construction, $\widetilde{G}$ contains the edge $u'v'$. Let $F$ be the set of such edges for every edge $uv\in E$. Note that, $\widetilde{G}/F = u(G)$ and that $\widetilde{G}/F$ is obtained by adding true twins and contracting edges. Since $\mathcal{C}$ is replicable and induced minor closed, $u(G)$ belongs to $\mathcal{C}$.
\end{proof}

\begin{obs}
The conditions given in Theorem~\ref{the:utterstable} are not necessary.    
\end{obs}
\noindent\begin{proof}
Consider the sequence $G_{i+1} = u(G_i)$, for $i \in \mathbb{Z}^+$, where $G_0 = C_4$. Let $\mathcal{C}$ be the set of all $G_i$ completed by induction with the graphs whose utter graph is in $\mathcal{C}$. Then, for any graph $G$ and its utter graph $u(G)$, we have that $G \in \mathcal{C}$ if and only if $u(G) \in \mathcal{C}$.
Note that each $G_i$ contains a $C_4$ and by construction, if $G$ is in $\mathcal{C}$, then there exists $i \in \mathbb{Z}^+$ such that, $G = G_i$ or there exists  $k \in \mathbb{Z}^+$ such that $u^k(G) = G_i$. Since $G$ is a clique if and only if $u(G)$ is a clique, no element of $\mathcal{C}$ is a clique, in particular $K_2 = (\{v_1,v_2 \}, \{ v_1v_2\})$ is not in $\mathcal{C}$ but is an induced subgraph of every element in $\mathcal{C}$. Hence $\mathcal{C}$ is not induced minor closed.
\end{proof} 

Since the class of contraction perfect graphs is replicable by  Lemma~\ref{lem:replication} and induced minor closed by definition, the following equivalence can be deduced from Theorem~\ref{the:utterstable}.
    
\begin{corollary}\label{the:equivcontractbothsides}
A graph is contraction perfect if and only if its utter graph is contraction perfect. 
\end{corollary}

Now, we show that contraction perfect graphs are precisely those for which the utter graph is perfect. This means that every perfect utter graph is also contraction perfect.

\begin{corollary}
\label{the:matching}
A graph is contraction perfect if and only if its utter graph is perfect. 
\end{corollary}

\noindent\begin{proof}
Let $G=(V,E)$ be a graph.

\noindent $(\Leftarrow)$ By definition of utter graphs, $G$ and $G/e$ for $e\in E $ are induced subgraphs of $u(G)$. Hence, all those graphs are perfect since $u(G)$ is.
By Theorem~\ref{the:contractequiv}, this implies that $G$ is contraction perfect.

\noindent $(\Rightarrow)$ If $G$ is contraction perfect, then so is $u(G)$ by Theorem~\ref{the:equivcontractbothsides}. By definition of contraction perfect graphs, this implies that $u(G)$ is perfect.
\end{proof}

Given a cost function $c$ on the vertices of $G$, let $\widetilde{c}_u = c_u$, for all $u\in V$, and $\widetilde{c}_{uv} = c_u + c_v$, for all $uv \in E$. Then finding a maximum $c$-weighted co-2-plex in $G$  is equivalent to finding a maximum $\widetilde{c}$-weighted stable set in $u(G)$, and then can be done in polynomial time ~\cite{GROTSCHEL1984325}. This gives the following corollary by Theorem~\ref{the:matching}.

\begin{corollary}\label{polynomial timetime}
    Finding a maximum weight co-2-plex in a contraction perfect graph is solvable in polynomial time.
\end{corollary}
%


%
 Therefore when $u(G)$ is perfect, computing a minimum covering by stable sets in $u(G)$ can be done in polynomial time. Note that a minimum cover of a graph by stable sets is exactly a graph coloring.
A further question is whether a graph coloring of $u(G)$ helps us to find a minimum covering of $G$ by co-2-plexes.
However, coloring the utter graph yields several coverings of $G$. Two opposite examples can be obtained as follows: a covering using only the stable sets of $u(G)[V]$, which are stable sets of $G$; and a covering using the stables of $u(G)[E]$ which are sets of edges of $G$ pairwise nonadjacent by contraction.

In Figure~\ref{color}, a covering of $G$ by co-$2$-plexes involves two co-$2$-plexes, whereas a coloring of $u(G)$ is composed of four stables sets. Hence the two problems are not equivalent.
An open question is whether a minimum co-2-plex covering of $G$ can be deduced from a given coloring of $u(G)$. Indeed, in Figure~\ref{color}, both examples given in the previous paragraphs yield  minimal coverings of $G$ $\{\{1,3\}, \{2\}\}$ and $\{\{12\}, \{23\}\}$.


\subsection{Subclasses of perfect graphs}\label{ss22}

This section shows how Theorem~\ref{the:utterstable} applies to other graph classes.
A \textit{split graph} is a graph whose vertex set can be partitioned into a stable set and a clique. A split graph is equivalently, a graph that does not contain $C_4$ or its complement.
A set of intervals is \emph{nested} if every couple of interval either has an empty intersection or one is contained in the other. A trivially perfect graph is an interval graph built from a nested interval set.
Equivalently, a \textit{trivially perfect graph} is a graph with no hole or path of size 4. 
For $k\ge 3$, a $k$-hole-free graph is a graph whose holes are all of size at most $k$. The well-known \textit{chordal graph} are the $3$-hole-free graphs.

Trivially perfect graphs are interval graphs, split and interval graphs are chordal. All these graph classes are contraction perfect, which is not the case for $k$-hole-free graphs with $k\ge 4$.

A graph class $\mathcal{C}$ is said {\it stable by utter graphs} if: $G\in\mathcal{C}$ if and only if $u(G)\in\mathcal{C}$.


\medskip

\begin{corollary}
The split (resp. trivially perfect, interval, chordal, $k$-hole-free) graphs are stable by utter graphs.
However, perfect graphs are not stable by utter graphs.
\end{corollary}
\noindent\begin{proof}
    One can show that all these graph classes are induced minor closed and replicable. We only detail the case of  interval graphs.
    
    Every induced subgraph $G_\mathcal{I}[W]$ of an interval graph $G_\mathcal{I}$ is the interval graph built from the interval set containing intervals of $\mathcal{I}$ associated with $W$.
    Contracting an edge $uv$ in $G_\mathcal{I}$ corresponds to replacing the intervals associated with $u$ and $v$ in $\mathcal{I}$ by the union of both intervals. 
    Adding a true twin to a vertex $v$ in an interval graph $G_\mathcal{I}$ corresponds to duplicating the interval associated with $v$ in $\mathcal{I}$.     
    Then, from Theorem~\ref{the:utterstable}, interval graphs are stable by utter graphs.
   
    Note that, for the $k$-hole-free graphs, the induced minor closeness follows from their definition and Lemma~\ref{lem:holegt5}.
    A counter example for a perfect graph to be not stable by utter graphs is $C_6$ which is perfect, but whose utter graph contains many $C_5$.
\end{proof}

Since the complement of perfect graphs are perfect, a natural question is to investigate the complement of contraction perfect graphs.
However, $\overline{C_6}$ is an antihole of size 6, hence contraction perfect, but $C_6$ is not.

\begin{obs} 
The complement of a contraction perfect graph is not necessarily contraction perfect. 
\end{obs}


\section*{Conclusion}

In this paper, we approach edge contraction in perfect graphs by characterizing the contraction perfect graphs in several manners: by contracting either a single edge or an edge set, by forbidden induced subgraphs and by the utter graph perfection.
As a byproduct, the recognition of contraction perfect graphs is polynomial.
Moreover, using utter graphs, the maximum weight co-2-plex problem becomes solvable in polynomial time in contraction perfect graphs. Finally, we focus on particular subclasses of perfect graphs.

For $k=2$, it remains an intriguing question whether a coloring of the utter graph can be used to get a minimal covering of the starting graph with co-$2$-plexes. 
When $k>2$, another interesting question is whether some extension of utter graphs could capture co-$k$-plexes as familiar combinatorial objects, like stable sets in utter graphs do for co-2-plexes.
From a polyhedral point of view, the equivalence between co-2-plexes of a graph and stable sets of its utter graph will give extended formulations for the co-2-plexe polytope from the stable set polytope.

Contraction perfect graphs being a new subclass of perfect graphs, one may also be interested in combinatorial algorithms to find a maximum clique/stable set or a minimum coloring on such graphs. Since, unlike for perfect graphs,  the complement of a contraction perfect graph is not necessarily contraction perfect, finding a maximum clique or a stable set may lead to distinct studies.

\section*{Acknowledgments}

We thank Daria Pchelina for her contribution to Lemma~\ref{lem:patates}.

\bibliographystyle{acm}
\bibliography{bibtex.bib}
\end{document}